\documentclass[11pt]{article}
\usepackage{mathrsfs}
\usepackage{amsfonts}
\usepackage{psfrag,epsfig}
\usepackage{latexsym,amsmath,amssymb}
\usepackage{exscale,relsize}
\usepackage{siunitx}
\usepackage{graphicx}
\usepackage{epstopdf}
\usepackage{indentfirst,cite}
\numberwithin{equation}{section}

\begin{document}
\title{{Non-homogeneous Schr\"{o}dinger systems with sign-changing and general nonlinearities: Infinitely many solutions}
\thanks{Research supported by Taishan Scholar Foundation for Young Experts of Shandong Province
(No. tsqn202306223) and Natural Science Foundation of Shandong Province (No. ZR2024MA009).}
\author {{Guanwei Chen$^{*}$}\\
\small  \it School of Mathematical Sciences, University of Jinan, Jinan 250022, Shandong Province, China\\
\small E-mail: guanweic@163.com}}
\date{}
\maketitle   \baselineskip 23pt

\date{}
\maketitle \baselineskip=12pt \noindent {\small {\bf Abstract}: In this paper, we study the non-homogeneous nonlinear Schr\"{o}dinger system
$$\left\{
\begin{array}{ll}
-\triangle u_j+V_j(x) u_j=g_j(x,u_1,\cdots,u_m)+h_j(x),& x\in \Omega,\\
\\
u_j:=u_j(x)=0,& x\in  \partial\Omega,\\
\\
j=1,2,\cdots,m,
\end{array}\right.
$$
where $\Omega\subset\mathbb{R}^{N}$ ($N\ge2$) is a bounded smooth domain, $(g_1,\cdots,g_m)$ is the gradient of $G(x,U)\in C^1(\Omega\times\mathbb{R}^m,\mathbb{R})$, $G(x,U)$ may be sign-changing, and it is super-quadratic or asymptotically-quadratic as $|U|\to\infty$. We obtain infinitely many solutions by using variational methods and perturbation methods, and we provide several typical examples to illustrate the main results. The {\bf main novelties} are as follows. (1) The nonlinearity $G$ may be sign-changing. (2) The nonlinearity $G$ is not only general, but also super-quadratic or asymptotically-quadratic at infinity and zero. (3) The nonlinearity $G$ is power-type or non-power-type. (4) We not only construct some new conditions, but also apply some conditions used in homogeneous problems to the study of non-homogeneous systems for the first time. The {\bf main difficulties} come from the following three aspects. (1) The proof of boundedness for $(PS)$ sequence of approximate functionals. (2) The detailed analysis of the asymptotic behaviors of approximate functionals. (3) The estimate of the upper and lower bounds for the minimax value sequence $\{c_k\}$ of the even function.



\baselineskip=18pt
\linespread{2.0}

\section*{1. Introduction and main results}
\setcounter{equation}{0} In this paper, we study the following non-homogeneous nonlinear Schr\"{o}dinger system
$$\left\{
\begin{array}{ll}
-\triangle u_j+V_j(x) u_j=g_j(x,u_1,\cdots,u_m)+h_j(x),& x\in \Omega,\\
\\
u_j:=u_j(x)=0,& x\in \partial\Omega,\\
\\
j=1,2,\cdots,m,
\end{array}\right.
\eqno(1.1)
$$
where $\triangle$ denotes the Laplace operator, $h_j(x)\not\equiv0$ are the non-homogeneous terms, $\Omega\subset\mathbb{R}^{N}$ ($N\ge2$) is a bounded smooth domain, and $(g_1,\cdots,g_m)$ is the gradient of $G(x,U)\in C^1(\Omega\times\mathbb{R}^m,\mathbb{R})$, i.e.,
$$
\nabla G(x,U)=\left(\frac{\partial G}{\partial u_1},\cdots,\frac{\partial G}{\partial u_m}\right)=(g_1,\cdots,g_m),\quad U=(u_1,\cdots,u_m)\in \mathbb{R}^{m}.
$$

The system (1.1) comes from the study the standing waves for the time-dependent  Schr\"{o}dinger system
$$\left\{
\begin{array}{ll}
i\frac{\partial\phi_j}{\partial t}=-\triangle\phi_j+\xi_j(x)\phi_j-f_j(x,|\phi_1|,\cdots,|\phi_m|)\phi_j, & t>0,\quad x\in \Omega,\\
\\
\phi_j:=\phi_j(t,x)=0, &  t>0,\quad x\in \partial\Omega,\\
\\
j=1,2,\cdots,m,
\end{array}\right.
\eqno(1.2)
$$
where $i$ is the imaginary unit.
By the definitions of standing waves ($\phi_j=e^{-i\omega_jt}u_j(x)$, $j=1,2,\ldots,m$), we get (1.2) becomes to
$$\left\{
\begin{array}{ll}
-\triangle u_j+\left(\xi_j(x)-\omega_j\right) u_j=f_j(x,|u_1|,\cdots,|u_m|)u_j,& x\in \Omega,\\
\\
u_j(x)=0, & x\in \partial\Omega,\\
\\
j=1,2,\cdots,m,
\end{array}\right.
\eqno(1.3)
$$
which is a special case of (1.1) with
$g_j(x,u_1,\cdots,u_m)\equiv f_j(x,|u_1|,\cdots,|u_m|)u_j$, $V_j\equiv \xi_j(x)-\omega_j$ and $h_j\equiv0$.

Nonlinear Schr\"{o}dinger systems can naturally model many physics-related problems, such as Bose-Einstein condensates and nonlinear optics \cite{Akhmediev,Deconinck,Hall,Ruegg,Timmermans}. Therefore, in the past decades, the existence and multiplicity of solutions of homogeneous Schr\"{o}dinger systems have been studied by many authors \cite{Costa-1,Bartsch1,Chang,Chen-0,Chen,Chen172,Chen2024,Ma-6,Ma-7,Pomponio,Zou-14,Zou-15,Zou-16}.

We note that the existing {\bf multiple results} of the homogeneous nonlinear Schr\"{o}dinger system are all based on the {\bf even condition} $I(U)=I(-U)$, where $I$ is the related energy functional. {\bf However}, the energy functional of (1.1) is not even (i.e., $\Phi(U)\neq\Phi(-U)$) because of the existence of the non-homogeneous terms $h_j(x)$, $j=1,\cdots,m$. Thus, the corresponding ideas based on the even functional {\bf are failed}. Inspired by \cite{Bolle-1,Bolle}, we study a continuous path of the related functionals, which starts from a symmetric (even) functional $\Phi_0$, and then prove a preservation result for the minmax critical levels to obtain critical points for the true functional $\Phi_1$ of (1.1).

Let $(\cdot,\cdot)$ and $|\cdot|$ denote respectively the usual inner product and the norm in $\mathbb{R}^{m}$. Let $\tau_t$ denote the Sobolev embedding constants such that
$$
\|U\|_t\le \tau_t\|U\|,\quad \forall U\in[H_{0}^{1}(\Omega)]^m, \quad \forall t\in[1,2^*),
\eqno(1.4)
$$
where $\|\cdot\|_t$ and $\|\cdot\|$ denote the norms in $[L^t(\Omega)]^m:=L^t(\Omega)\times\cdots\times L^t(\Omega)$ and $[H_{0}^{1}(\Omega)]^m:=H_{0}^{1}(\Omega)\times\cdots\times H_{0}^{1}(\Omega)$.\\

$\bigstar$ First, we study the case of {\it super-quadratic} growth at infinity. Assume that\\
$\bf{(V_1)}\quad$There are $\underline{V},\overline{V}>0$ such that
$\underline{V}\le V_j(x)\le \overline{V}$, $\forall x \in\Omega$, where $j=1,\cdots,m.$\\
$\bf{(G_{1})}\quad$$G(x,U)\in C^1(\Omega\times\mathbb{R}^m,\mathbb{R})$ and there exist $\varepsilon_0,\widetilde{\varepsilon_0}>0$ such that
$$
\frac{|\nabla G(x,U)|}{|U|}\le \varepsilon_0<\min\left\{1/\tau_2^2,\underline{V}\right\},\quad \forall|U|\le\widetilde{\varepsilon_0},\quad \forall x\in\Omega.
$$
$\bf{(G_{2})}\quad$$\widetilde{G}(x,U):=\frac{1}{2}(\nabla G(x,U),U)-G(x,U)>0$, $\forall(x,U)\in\Omega\times(\mathbb{R}^m\setminus\{0\})$.\\
$\bf{(G_{3})}\quad$There are $a_1,a_2>0$ and $\mu\in(1,2^\ast)$  ($2^\ast:=\frac{2N}{N-2}$ if $N\ge3$ and $2^\ast:=+\infty$ if $N=2$) such that
$$
\widetilde{G}(x,U) \ge a_1|U|^\mu,\quad \forall |U|\ge a_2,\quad \forall x \in\Omega.
$$
$\bf{(G_{4})}\quad$$\mathlarger{\lim_{|U|\rightarrow\infty}}\frac{G(x,U)}{|U|^{2}}=+\infty$ uniformly in $x\in\Omega$.\\
$\bf{(G_{5})}\quad$There are $a_3,a_4>0$ and $p\in(2,2^\ast)$ such that
$$
|\nabla G(x,U)|\le a_3|U|^{p-1}, \quad\forall |U|\ge a_4, \quad\forall x\in\Omega.
$$
$\bf{(G_{6})}\quad$There exist $r_0,r_1>0$ and $\gamma_1,\gamma_2\in(0,1)$ such that
$$
(\nabla G(x,U),U)\leq r_0\left(\widetilde{G}(x,U)\right)^{\gamma_1}|U|^{\gamma_2}, \quad\forall|U|\geq r_1, \quad \forall x\in\Omega.
$$
$\bf{(H_1)}\quad$$h_j\in L^{\mu'}(\Omega)\setminus\{0\}$, where $\mu^{-1}+ (\mu')^{-1}=1$ and  $\mu\in(1,2^\ast)$,  $j=1,\cdots,m$.

Then we have the following results.
\\
$\bullet$ \textbf{Theorem 1.1.} {\it Assume $(V_1)$, $(H_{1})$ and $(G_{1})$-$(G_{6})$ hold.

$(1)$ If $N=2$ and the parameters $\mu$ and $p$ in $(G_{3})$ and $(G_{5})$ satisfy
$$
\mu>1+\frac{\varepsilon+1}{\varepsilon+2}>\frac{3}{2} \quad \mbox{and}\quad
2+\frac{1}{\varepsilon+1}\le p<\frac{2\mu(\varepsilon+1)}{\mu\varepsilon+1}, \qquad \forall\varepsilon>0,
\eqno(1.5)
$$
then the non-homogeneous Schr\"{o}dinger system $(1.1)$ has infinitely many solutions.

$(2)$ If $N\ge3$ and the parameters $\mu$ and $p$ in $(G_{3})$ and $(G_{5})$ satisfy
$$
\mu>1+\frac{N}{N+2}\ge\frac{8}{5}\quad \mbox{and}\quad 2+\frac{2}{N}\le p<\frac{2\mu N}{\mu(N-2)+2},
\eqno(1.6)
$$
then the non-homogeneous Schr\"{o}dinger system $(1.1)$ has infinitely many solutions.}
\\
$\bullet$ \textbf{Theorem 1.2.} {\it Assume that $(V_1)$, $(H_{1})$, $(G_{1})$-$(G_{5})$ and the following condition hold\\
$\bf{(G'_{6})}\quad$There exist $\mu_0>2$ and $\kappa_0>0$ such that
$$
(\nabla G(x,U),U)\ge \mu_0 G(x,U)-\kappa_0 |U|^{2},\quad\forall(x,U)\in\Omega\times\mathbb{R}^m.
$$
Then we also have the same results $(1)$ and $(2)$ in Theorem $1.1$.}
\\
$\bullet$ \textbf{Corollary 1.1.} {\it Assume that $(V_1)$, $(H_{1})$, $(G_{1})$, $(G_{5})$ and the following condition hold\\
$\bf{(G''_{6})}\quad$There is $\mu_0>2$ such that
$$
(\nabla G(x,U),U)\ge \mu_0 G(x,U)>0,\quad \forall (x,U)\in\Omega\times(\mathbb{R}^m\setminus\{0\}).
$$
Then we also have the same results $(1)$ and $(2)$ in Theorem $1.1$.}
\\
$\bullet$ \textbf{Remark 1.1.} Note that $(G''_{6})$ implies the conditions $(G_{2})$-$(G_{4})$ and $(G'_{6})$ hold; thus, Corollary 1.1 can be proved by Theorem 1.2. The condition $(G_{1})$ implies $G$ may be super-quadratic or asymptotically-quadratic at zero, i.e.,
$$
\mathlarger{\lim_{|U|\rightarrow0}}\frac{G(x,U)}{|U|^2}=l(x),\quad l(x)\equiv0 \ \mbox{or} \ l(x)\in L^{\infty}(\Omega)\setminus\{0\},
$$
see the following functions in Example 1.1. The roles of the conditions $(G_{6})$ and $(G'_{6})$ are mainly used to prove the boundedness of $(PS)$ sequences of the approximate functionals $\Phi_\theta$, see the following Lemma 2.2.
\\
$\bullet$ \textbf{Example 1.1.} To illustrate Theorems 1.1-1.2 and Corollary 1.1, we give four examples.

{\it Ex1 (super-quadratic at infinity and zero, non-power-type).}
$$
G(x,U)=\eta_1(x)\left[(\delta-2)|U|^{\delta-\varepsilon}\sin^{2}\frac{|U|^{\varepsilon}}{\varepsilon}
+|U|^{\delta}\right], \quad \delta>2, \ 0<\varepsilon<\delta-2.
$$

{\it Ex2 (super-quadratic at infinity and zero, non-power-type).}
$$
G(x,U)=\eta_2(x)|U|^2\ln\left(1+|U|^2\right).
$$

{\it Ex3 (super-quadratic at infinity, asymptotically-quadratic at zero, non-power-type, sign-changing).}
$$
G(x,U)=\eta_3(x)|U|^{2}\left(|U|^{2}-2\cos |U|\right), \quad 4\left|\sup_{x\in\Omega}\eta_3(x)\right|<\underline{V}.
$$

{\it Ex4 (super-quadratic at infinity and zero, power-type).}
$$
G(x,U)=\eta_4(x)|U|^p, \quad p>2.
$$
Here, $(x,U)\in\Omega\times\emph{R}^m$ and $0<\mathlarger{\inf_{x\in\Omega}}\eta_j(x)\le\mathlarger{\sup_{x\in\Omega}}\eta_j(x)<+\infty$, $j=1,2,3,4$.
By calculation, we can prove that the functions in (Ex1)-(Ex4) all satisfy our conditions $(G_1)$-$(G_6)$, the functions in (Ex3) and (Ex4) also satisfy our condition $(G'_6)$, and the function in (Ex4) also satisfies our condition $(G''_6)$, but the functions in (Ex1)-(Ex3) all don't satisfy our condition $(G''_6)$. Clearly, the function $G(x,U)$ in (Ex3) is sign-changing, that is, $G(x,U)<0$ if $|U|$ is small and $G(x,U)>0$ if $|U|$ is large.\\

$\bigstar$ Second, we study the case of {\it asymptotically-quadratic} growth at infinity. Assume that\\
$\bf{(G_{7})}\quad$There is $\eta(x)\in L^{\infty}(\Omega)$ such that $\mathlarger{\inf_{x\in\Omega}}\eta(x)>0$ is large enough and
$$
\lim_{|U|\rightarrow\infty}\frac{|\nabla G(x,U)|}{|U|}
=2\lim_{|U|\rightarrow\infty}\frac{G(x,U)}{|U|^2}
=\eta(x) \quad \mbox{uniformly in} \ x\in\Omega.
$$
$\bf{(G_{8})}\quad$There is $\mu\in(1,2^\ast)$ such that
$\mathlarger{\lim_{|U| \to \infty}}\frac{\widetilde{G}(x,U)}{|U|^{\mu}}=+\infty$ uniformly in $x\in\Omega.$

Then we have the following results.
\\
$\bullet$ \textbf{Theorem 1.3.} {\it Assume that $(V_1)$, $(G_{1})$, $(G_{2})$, $(G_{7})$, $(G_{8})$ and $(H_{1})$ with $\mu$ is given in $(G_{8})$ hold.

$(1)$ If $N=2$ and the parameter $\mu$ in $(G_{8})$ satisfies
$$
\mu>1+\frac{\varepsilon+1}{\varepsilon+2}>\frac{3}{2}, \quad \forall\varepsilon>0,
\eqno(1.7)
$$
then the non-homogeneous Schr\"{o}dinger system $(1.1)$ has infinitely many solutions.

$(2)$ If $N\ge3$ and the parameter $\mu$ in $(G_{8})$ satisfies
$$
\mu>1+\frac{N}{N+2}\ge\frac{8}{5},
\eqno(1.8)
$$
then the non-homogeneous Schr\"{o}dinger system $(1.1)$ has infinitely many solutions.}
\\
$\bullet$ \textbf{Example 1.2.} To illustrate Theorems 1.3, we give the following example.

{\it Ex5 (asymptotically-quadratic at infinity, super-quadratic at zero, non-power-type).}
$$
G(x,U)=\eta(x)\left[\tau-\frac{\tau}{\ln(e+|U|)}\right]|U|^{2},\quad \forall(x,U)\in\Omega\times \emph{R}^{m},
$$
where $0\leq \eta(x)\in L^{\infty}(\Omega)\setminus\{0\}$ and $\tau$ is a positive constant. By calculation, we can prove that it satisfies our conditions $(G_{1})$, $(G_{2})$, $(G_{7})$ and $(G_{8})$.
\\
$\bullet$ \textbf{Remark 1.2.} Because of the existence of non-homogeneous terms and general nonlinearities (especially  non-power-type), it is difficult to prove the existence of infinitely many solutions of the system (1.1). Specifically, the {\bf first difficulty} is the proof of boundedness for $(PS)$ sequence of approximate functionals (see Lemma 2.2); the {\bf second difficulty} is a detailed analysis of the asymptotic behaviors of approximate functionals (see Lemma 2.3); the {\bf third difficulty} is the estimate of lower bound for the minimax value sequence $\{c_k\}$ of the even function (see Section 3).
\\
$\bullet$ \textbf{Remark 1.3.} The {\bf main novelties} of this paper are as follows.
(1) the nonlinearity $G$ may be sign-changing (see (Ex3)).
(2) the nonlinearity $G$ is not only general, but also super-quadratic or asymptotically-quadratic at infinity and zero.
(3) the nonlinearity $G$ includes two cases: power-type (see (Ex4)) and non-power-type (see (Ex1)-(Ex3) and (Ex5)).
(4) although the conditions $(G'_{6})$ and $(G_{8})$ have been applied to the study of homogeneous problems, we apply them to the problem of proving multiple solutions of non-homogeneous systems for the first time, besides, we construct a new condition $(G_{6})$ to prove the boundedness of $(PS)$ sequence of approximate functionals in the non-homogeneous systems. In short, the conditions and methods used in this paper will provide useful ideas for the study of related non-homogeneous equations and systems.

The rest of the paper is organized as follows. In Section 2, we
establish the variational framework of the system (1.1), and
give some preliminary lemmas. In Section 3, we give the detailed proofs of our
main results.

\section * {2. Variational framework and preliminary lemmas}

We first give the Bolle's Perturbation theorem.
Let $W$ be a Hilbert space equipped with the norm $\|\cdot\|$. Let $W=W_- \oplus  W_+$ and $\dim(W_-) < \infty$, and let $\{e_n\}_{n\ge1}$ be an orthonormal basis of $W_+$.
Let
$$
W_0:=W_-, \quad W_{k+1}:=W_k \oplus \mathbb{R}e_{k+1}, \quad k=0,1,2,\cdots,
$$
then $\{W_k\}$ is an increasing sequence of finite dimensional subspaces of $W$. We consider the family of functionals $\Phi_{\theta}(u): \ [0,1]\times W\to \mathbb{R}$, which are assumed that they are $C^1$ according to \cite{Clapp}. \\
\\
\textbf{Lemma 2.1} (see \cite{Bolle-1,Bolle}). {\it Let $\Phi'_{\theta}(u)=\frac{\partial\Phi_{\theta}(u)}{\partial u}$. Assume $\Phi_{\theta}(u)$  satisfies the following conditions.\\
$(A1)$ $\{(\theta_n,u_n)\}\subset[0,1]\times W$ has a convergent subsequence if
$|\Phi_{\theta_n}(u_n)|<\infty$ and $\mathlarger{\lim_{n\to\infty}}\Phi'_{\theta_n}(u_n)= 0.$\\
$(A2)$ For any $\varepsilon>0$ there exists $C_\varepsilon> 0$ such that
$$
|\Phi_{\theta}(u)|\le\varepsilon\Rightarrow\left|\frac{\partial\Phi_{\theta}(u)}{\partial \theta}\right|\le C_\varepsilon(\|\Phi'_{\theta}(u)\|+1)(\|u\|+1), \quad \forall(\theta,u)\in[0,1]\times W.
$$
$(A3)$ There exist two continuous maps $\eta_1, \eta_2 : [0, 1]\times \mathbb{R}\to \mathbb{R}$ $($they are all Lipschitz continuous
with respect to the second variable$)$ such that
$$
\eta_1(\theta,\cdot)\le\eta_2(\theta,\cdot), \quad
\Phi'_{\theta}(u)=0 \Rightarrow \eta_1(\theta,\Phi_{\theta}(u))\le \frac{\partial\Phi_{\theta}(u)}{\partial \theta} \le \eta_2(\theta,\Phi_{\theta}(u)), \quad \forall(\theta,u)\in[0,1]\times W.
$$
$(A4)$ $\Phi_{0}(-u)=\Phi_{0}(u)$ and $\mathlarger{\lim_{u\in H,\ \|u\|\to\infty}\sup_{\theta\in[0,1]}}\Phi_{\theta}(u)=-\infty$, $\forall H\subset W$ with $\dim H<\infty.$

Then for all $k\in\mathbb{N}$ we have\\
$(i)$ either $\Phi_1$ has a critical level $\widetilde{c_k}$ with  $\psi_2(1,c_k)<\psi_1(1,c_{k+1})\le\widetilde{c_k}$,\\
$(ii)$ or $c_{k+1}-c_{k}\le C(\overline{\eta}_1(c_{k+1})+\overline{\eta}_2(c_{k})+1)$, where $c_k$, $\overline{\eta}_i$ and $\psi_i$ are defined as follows:
$$
c_k=\inf_{\gamma\in \Gamma_k}\sup_{u\in W_k}\Phi_0(\gamma(u)), \quad \Gamma_k=\{\gamma\in C(W_k,W):\ \gamma \ is \ odd, \ \gamma(u)=u\in W \ for \ \|u\| \ large\},
$$
$$
\overline{\eta}_1(s)=\sup_{\theta\in[0,1]}|\eta_1(\theta,s)|, \quad \overline{\eta}_2(s)=\sup_{\theta\in[0,1]}|\eta_2(\theta,s)|,
$$
and $\psi_j : [0,1]\times \mathbb{R},\mathbb{R}$ are the solutions
of the problem
$$
\left\{
\begin{array}{ll}
\frac{\partial\psi_j(\theta,s)}{\partial\theta}=\eta_j(\theta,\psi_j(\theta,s)),\\
\psi_j(0,s)=s,
\end{array}\right.
\eqno(2.0)
$$
where $\psi_j(\theta,\cdot)$ is continuous and $\psi_1(\theta,\cdot)\le\psi_2(\theta,\cdot)$, $j=1,2$.}\\

In what follows, we always assume the conditions in Theorems 1.1-1.3 always hold. Let $C_j$ ($j=0,1,\cdots$) denote different positive constant.

Let $\|\cdot\|_{L^t}$ and $\|\cdot\|_{t}$ denote respectively
the norms of $L^{t}:=L^{t}(\Omega)$ and $[L^{t}]^m:=L^{t}\times\cdots\times L^{t}$, where
$$
\|U\|_{t}=\left(\int_\Omega|U|^t dx\right)^\frac{1}{t}=\left(\int_\Omega(|u_1|^2+\cdots+|u_m|^2)^{\frac{t}{2}} dx\right)^\frac{1}{t}, \quad 1\le t<2^*,
$$
and
$$
\|U\|_{\infty}=\sup_{x\in\Omega}|U|=\sup_{x\in\Omega}\left(|u_1|^2+\cdots+|u_m|^2\right)^{\frac{1}{2}},  \quad U=(u_1,\cdots,u_m).
$$
Let $H_{0}^{1}(\Omega)$ be the usual Sobolev space and
$$
W:=[H_{0}^{1}(\Omega)]^m=H_{0}^{1}(\Omega)\times\cdots\times H_{0}^{1}(\Omega).
$$
By $(V_1)$, we can let $\|\cdot\|_{\overline{H}_i}$ be the norm of $H_{0}^{1}(\Omega)$ generated by the inner product
$$
(u,v)_{j}=\int_{\Omega}\left(\nabla u\nabla v+V_j(x)uv\right)dx, \quad  \forall u,v\in
H_{0}^{1}(\Omega), \quad j=1,\cdots,m.
$$
Then the induced inner product and norm on $W$ are given
respectively by
$$
\langle U,V\rangle=(u_1,v_1)_{1}+\cdots+(u_m,v_m)_{m}
$$
and
$$
\|U\|=\langle U,U\rangle^{1/2}, \quad \forall U,  V\in W,
$$
where $U=(u_{1},\cdots,u_{m})$ and $V=(v_{1},\cdots,v_{m})$.
We mention that the norms in (1.4) are the same as the norms $\|\cdot\|_t$ and $\|\cdot\|$ given above.
By the above definitions, we have
$$
\|u_j\|_{\overline{H}_j}\le\|U\|, \quad \|u_j\|_{L^t}\le\|U\|_{t}, \quad \forall t\in[1,\infty], \ \forall U=(u_{1},\cdots,u_{m})\in W, \quad   j=1,\cdots,m.
\eqno(2.1)
$$

Next, we consider the family of approximate functionals of the system (1.1):
$$\begin{array}{rcl}
\Phi_\theta(U)&=&\frac{1}{2}\|U\|^2-\mathlarger{\int_{\Omega}}G(x,U)dx
-\theta\mathlarger{\int_{\Omega}}(H,U)dx\\
\\
&=&\Phi_0(U)-\theta\mathlarger{\int_{\Omega}\sum_{j=1}^{m}}h_j(x)u_jdx,\quad \forall \theta\in[0,1], \quad \forall U \in W,
\end{array}\eqno(2.2)
$$
where $U=(u_{1},\cdots,u_{m})$, $H=(h_{1},\cdots,h_{m})$ and
$$
\Phi_0:=\frac{1}{2}\|U\|^2-\int_{\Omega}G(x,U)dx.
$$
Obviously, $\Phi_0$ is an even functional, and $\Phi_1$ is the ``ture" energy functional of (1.1).
From \cite{Mawhin}, we know
$\Phi_\theta \in C^{1}(W,\mathbb{R})$ and the derivatives are given by
$$\begin{array}{rcl}
\Phi'_\theta(U)V&=&\Phi'_\theta(u_{1},\cdots,u_{m})(v_{1},\cdots,v_{m})\\
\\
&=&\langle U,V\rangle-\mathlarger{\int_{\Omega}}(\nabla G(x,U),V)dx
-\theta\mathlarger{\int_{\Omega}\sum_{j=1}^{m}}h_j(x)v_jdx
\end{array}
\eqno(2.3)
$$
and
$$
\frac{\partial\Phi_{\theta}(U)}{\partial \theta}=-\int_{\Omega}\sum_{j=1}^{m}h_j(x)u_jdx,\quad
\forall U,V\in W, \ \ \forall\theta\in[0,1].
$$
Therefore, the (weak) solutions of the system (1.1) are the critical points of the ``true" functional $\Phi_1$.\\

To apply Lemma 2.1, we need to verify the related conditions $(A1)$-$(A4)$ in Lemma 2.1 with our definitions, i.e., the following Lemmas 2.2 and 2.3.\\
\\
\textbf{Lemma 2.2.} {\it $\Phi_{\theta}(U)$ satisfies the condition $(A1)$. }\\
\textbf{Proof.} Assume $\{(\theta_n,U^n)\}\subset[0,1]\times W$ satisfies
$$
|\Phi_{\theta_n}(U^n)|<\infty, \quad \Phi'_{\theta_n}(U^n) \to 0, \ \ n\to\infty,
\eqno(2.4)
$$
where $U^n=(u_1^n,\cdots,u_m^n)$. We will prove the sequence $\{(\theta_n,U^n)\}$ has a convergent subsequence. First, we prove the sequence $\{U^n\}$ satisfies $\|U^n\|<\infty$. Second, we prove $\{U^n\}$ has a convergent subsequence. Therefore, $\{(\theta_n,U^n)\}$ has a convergent subsequence since $\{\theta_n\}\subset[0,1]$.

(1) Frist, we prove $\|U^n\|<\infty$. We divide the proof of (1) into three parts by the assumptions in Theorems 1.1-1.3.

\textbf{Part 1.} If the assumptions in Theorem 1.1 hold. Suppose that $\|U^n\|\to\infty$. Let
$$
V^n:=\frac{U^n}{\|U^n\|},
$$
then $\|V^n\|=1$. After passing to a subsequence, then there is $V\in W$ such that
$$
V^n\rightharpoonup V \ \mbox{in} \ W, \qquad
V^n(x)\to V(x) \ \mbox{a.e.} \ x\in\Omega.
\eqno(2.5)
$$

By (2.3), we have
$$
\Phi'_{\theta_n}(U^n)U^n=\|U^n\|^2-\int_{\Omega}(\nabla G(x,U^n),U^n)dx-\theta\int_{\Omega}\sum_{j=1}^{m}h_j(x)u^n_jdx,
\eqno(2.6)
$$
and
$$
\frac{\Phi'_{\theta_n}(U^n)}{\|U^n\|^2}=1- \frac{\mathlarger{\int_{\Omega}}(\nabla G(x,U^n),U^n)dx}{\|U^n\|^2}
-\frac{\theta\mathlarger{\int_{\Omega}\sum_{j=1}^{m}}h_j(x)u^n_jdx}{\|U^n\|^2}.
\eqno(2.7)
$$
By the H\"{o}lder's inequality, the Sobolev embedding theorem, $(H_1)$, (2.1) and $\|U^n\|\to\infty$, we have
$$
\lim_{n\to\infty}\frac{\left|\theta\mathlarger{\int_{\Omega}\sum_{j=1}^{m}}
h_j(x)u^n_jdx\right|}{\|U^n\|^2}
\le\lim_{n\to\infty}\frac{\theta\left(\mathlarger{\sum_{j=1,\cdots,m}}
\|h_j\|_{L^{{\mu'}}}\right)\|U\|_{\mu}}{\|U^n\|^2}=0.
\eqno(2.8)
$$
Hence, by $\|U^n\|\to\infty$, $\Phi'_{\theta_n}(U^n) \to 0$ in (2.4) and (2.7)-(2.8), we have
$$
\lim_{n\to\infty} \frac{\mathlarger{\int_{\Omega}}(\nabla G(x,U^n),U^n)dx}{\|U^n\|^2}=1.
\eqno(2.9)
$$

For all $b>a\ge0$, we let
$$
\Omega_{n}(a,b):=\left\{x\in\Omega:
 \ a\leq|U^n|<b\right\}
$$
for the given functions $U^n$. By (2.2)-(2.4), the H\"{o}lder's inequality, (2.1) and the Sobolev embedding theorem, for $n$ is large enough, we have
$$\begin{array}{rcl}
C_0&\geq&\Phi_{\theta_n}(U^n)-\frac{1}{2}\Phi'_{\theta_n}(U^n)U^n\\
\\
&=&\mathlarger{\int_{\Omega}}\widetilde{G}(x,U^n)dx
-\frac{\theta}{2}\mathlarger{\int_{\Omega}\sum_{j=1}^{m}}h_j(x)u^n_jdx\\
\\
&\ge&\mathlarger{\int_{\Omega}}\widetilde{G}(x,U^n)dx-\frac{1}{2} \left(\mathlarger{\sum_{j=1}^{m}}\|h_j\|_{L^{\mu'}}\right)\|U^n\|_{\mu}\\
\\
&\ge&\mathlarger{\int_{\Omega}}\widetilde{G}(x,U^n)dx-C_1\|U^n\|,
\end{array}
\eqno(2.10)
$$
thus for $n$ is large enough, we have
$$\begin{array}{rcl}
\mathlarger{\int_{\Omega}}\widetilde{G}(x,U^n)dx
&=&\mathlarger{\int_{x\in\Omega_{n}(0,a)}\widetilde{G}(x,U^n)dx
+\int_{x\in\Omega_{n}(a,b)}\widetilde{G}(x,U^n)dx +
\int_{x\in\Omega_{n}(b,+\infty)}\widetilde{G}(x,U^n)dx}\\
\\
&\le& C_2(1+\|U^n\|).
\end{array}
\eqno(2.11)
$$

By (1.4) and $\|V^n\|=1$, we have
$$
\|V^n\|^2_2\le \tau_2^2\|V^n\|^2=\tau_2^2.
$$
It follows from $(G_1)$ that
$$\begin{array}{rcl}
\frac{\mathlarger{\int_{x\in\Omega_{n}(0,\widetilde{\varepsilon_0})}}(\nabla G(x,U^n),U^n)dx}{\|U^n\|^2}
&\le&\mathlarger{\int_{x\in\Omega_{n}(0,\widetilde{\varepsilon_0})}}\frac{|\nabla G(x,U^n)|}{|U^n|}|V^n|^2dx\\
\\
&\leq&\varepsilon_0\|V^n\|^2_2 \\
\\
&\leq&\varepsilon_0\tau_2^2<1, \quad \forall n\in \emph{N}^+.
\end{array}
\eqno(2.12)
$$
By (1.4) and $\|V^n\|=1$, we have
$$
\|V^n\|_{\frac{1}{1-\gamma_1}}\le \tau_{\frac{1}{1-\gamma_1}},
\eqno(2.13)
$$
where $\gamma_1\in(0,1)$ is given in $(G_6)$, and $\frac{1}{1-\gamma_1}>1$.

Let
$$
\varepsilon\in\left(0,\frac{1-\varepsilon_0\tau_2^2}{4}\right).
$$
By $(G_6)$, (2.11), (2.13) and $\|U^n\|\to\infty$, we can take $b_\varepsilon> \widetilde{\varepsilon_0}$ such that
$$\begin{array}{rcl}
\frac{\mathlarger{\int_{x\in\Omega_{n}(b_\varepsilon,+\infty)}}(\nabla G(x,U^n),U^n)dx}{\|U^n\|^2}
&\leq&\mathlarger{\int_{x\in\Omega_{n}(b_\varepsilon,+\infty)}}
\frac{r_0\left(\widetilde{G}(x,U^n)\right)^{\gamma_1}|U^n|^{\gamma_2}|V^n|}{|U^n|\|U^n\|}dx\\
\\
&=&\mathlarger{\int_{x\in\Omega_{n}(b_\varepsilon,+\infty)}}
\frac{r_0\left(\widetilde{G}(x,U^n)\right)^{\gamma_1}|V^n|}{|U^n|^{1-\gamma_2}\|U^n\|}dx\\
\\
&\le&\frac{r_0}{b_\varepsilon^{1-\gamma_2}}
\frac{\mathlarger{\int_{x\in\Omega_{n}(b_\varepsilon,+\infty)}}
\left(\widetilde{G}(x,U^n)\right)^{\gamma_1}|V^n|dx}{\|U^n\|}\\
\\
&\leq&\frac{r_0}{b_\varepsilon^{1-\gamma_2}}
\frac{\left(\mathlarger{\int_{x\in\Omega_{n}(b_\varepsilon,+\infty)}}
\widetilde{G}(x,U^n)dx\right)^{\gamma_1}\|V^n\|_{\frac{1}{1-\gamma_1}}}{\|U^n\|}\\
\\
&\leq&\frac{r_0}{b_\varepsilon^{1-\gamma_2}}
\frac{C_3\left(1+\|U^{n}\|\right)^{\gamma_1}}{\|U^{n}\|} \qquad \left(C_3:=C_2^{\gamma_1}\tau_{\frac{1}{1-\gamma_1}}\right)\\
\\
&<&\varepsilon<\frac{1-\varepsilon_0\tau_2^2}{4},\quad \forall n\ge K_1,
\end{array}
\eqno(2.14)
$$
for some $K_1>0$. Let
$$
C^{b_\varepsilon}_{\widetilde{\varepsilon_0}}:=\inf\left\{\frac{\widetilde{G}(x,U)}{|U|^2}:\ x\in\Omega,
 \ U\in\mathbb{R}^{m}, \ \widetilde{\varepsilon_0}\leq|U|\le b_\varepsilon\right\}.
$$
Note that $\widetilde{G}\in C(\Omega\times\mathbb{R}^m,\mathbb{R})$ and $(G_2)$ imply that there exists $\widetilde{c_0}>0$ such that
$$
\widetilde{G}(x,U)\ge \widetilde{c_0}> 0 \quad \mbox{if} \ \widetilde{\varepsilon_0}\leq|U|\le b_\varepsilon,
$$
so we have $C^{b_\varepsilon}_{\widetilde{\varepsilon_0}}\ge \frac{\widetilde{c_0}}{|b_\varepsilon|^2}> 0$ and
$$
\widetilde{G}(x,U^n)\geq C^{b_\varepsilon}_{\widetilde{\varepsilon_0}}|U^n|^2,\quad\forall
x\in\Omega_{n}(\widetilde{\varepsilon_0},b_\varepsilon).
$$
It follows from (2.11) and $\|U^n\|\to\infty$ that
$$\begin{array}{rcl}
\mathlarger{\int_{x\in\Omega_{n}(\widetilde{\varepsilon_0},b_\varepsilon)}}|V^n|^2dx
&=&\frac{\mathlarger{\int_{x\in\Omega_{n}(\widetilde{\varepsilon_0},b_\varepsilon)}}|U^n|^2dx}{\|U^n\|^2}\\
\\
&\le&\frac{\mathlarger{\int_{x\in\Omega_{n}(\widetilde{\varepsilon_0},b_\varepsilon)}}\widetilde{G}(x,U^n)dx}
{C^{b_\varepsilon}_{\widetilde{\varepsilon_0}}\|U^n\|^2}\\
\\
&\le&\frac{C_2(1+\|U^n\|)}{C^{b_\varepsilon}_{\widetilde{\varepsilon_0}}\|U^n\|^2}\rightarrow0, \quad n\rightarrow\infty.
\end{array}
\eqno(2.15)
$$
By $(G_1)$ and $|\nabla G|\in C(\Omega\times\mathbb{R}^m,\mathbb{R})$, there is $\varrho=\varrho(\widetilde{\varepsilon_0},b_\varepsilon)>0$
independent of $n$ such that
$$
|\nabla G(x,U^n)|\leq\varrho
|U^n|,\quad \forall x\in\Omega_{n}(\widetilde{\varepsilon_0},b_\varepsilon).
$$
It follows from (2.15) that there is
$K_2>0$ such that
$$\begin{array}{rcl}
\frac{\mathlarger{\int_{x\in\Omega_{n}(\widetilde{\varepsilon_0},b_\varepsilon)}}(\nabla G(x,U^n),U^n)dx}{\|U^n\|^2}
&\le& \varrho\int_{x\in\Omega_{n}(\widetilde{\varepsilon_0},b_\varepsilon)}|V^n|^2dx\\
\\
&<& \varepsilon<\frac{1-\varepsilon_0\tau_2^2}{4},\quad \forall n\geq K_2.
\end{array}
\eqno(2.16)
$$
Hence, by (2.12), (2.14), (2.16) and $0<\varepsilon_0\tau_2^2<1$ by $(G_1)$, we have
$$
\lim_{n\to\infty} \frac{\mathlarger{\int_{x\in\Omega}}(\nabla G(x,U^n),U^n)dx}{\|U^n\|^2}\le\frac{1+\varepsilon_0\tau_2^2}{2}<1, \quad \forall n\geq\max\{K_1,K_2\},
$$
which contradicts with (2.9).

Therefore, $\{U^n\}$ is bounded in $W$, i.e., $\|U^n\|<\infty$.

\textbf{Part 2.} If the assumptions in Theorem 1.2 hold. If $\|U^n\|\rightarrow\infty$, we let
$$
V^n:=\frac{U^n}{\|U^n\|},
$$
then $\|V^n\|=1$ and (2.5) still holds.  By (2.1)-(2.4), $(G'_6)$, $(H_1)$, the H\"{o}lder's inequality and the Sobolev embedding theorem,  we have
$$
\begin{array}{rcl}
C_4&\ge&  \Phi_{\theta_n}(U^n)- \frac{1}{\mu_0}\Phi'_{\theta_n}(U^n)U^n\\
\\
&=& \left(\frac{1}{2}-\frac{1}{\mu_0}\right)\|U^n\|^2+\mathlarger{\int_{\Omega}}\left[\frac{1}{\mu_0}(\nabla G(x,U^n),U^n)-G(x,U^n)\right]dx\\
\\
&&-\theta_n\left(1-\frac{1}{\mu_0}\right)\mathlarger{\int_{\Omega}\sum_{j=1}^{m}}h_j(x)u^n_jdx\\
\\
&\ge& \left(\frac{1}{2}-\frac{1}{\mu_0}\right)\|U^n\|^2-\frac{\kappa_0}{\mu_0} \|U^n\|^{2}_2
-\theta_n\left(1-\frac{1}{\mu_0}\right)\left(\mathlarger{\sum_{j=1,\cdots,m}}
\|h_j\|_{L^{{\mu'}}}\right)\|U\|_{\mu}\\
\\
&\ge& \left(\frac{1}{2}-\frac{1}{\mu_0}\right)\|U^n\|^2-\frac{\kappa_0}{\mu_0} \|U^n\|^{2}_2-C_5\|U^n\|,  \quad \forall n\in \mathbb{N}^+,
\end{array}
$$
which together with $\|U^n\|\rightarrow\infty$ and $V^n=\frac{U^n}{\|U^n\|}$ implies
$$
\limsup_{n\to\infty}\|V^n\|^{2}_2\ge \frac{\mu_0-2}{2\kappa_0}>0.
$$
It follows from (2.5) and the Sobolev compact embedding theorem that
$$
V\neq0.
$$

Since $(G_4)$ or $(G_7)$ implies that there are $a_0,\zeta_0>0$ such that
$$
G(x,U)\ge\zeta_0|U|^2,\quad \forall|U|\ge a_0>0.
\eqno(2.17)
$$
By $\|U^n\|\rightarrow\infty$, $V^n=\frac{U^n}{\|U^n\|}$ and (2.5), we have
$$
\lim_{n\rightarrow\infty}|U^{n}|=+\infty,\quad
\forall x\in \widetilde{\Omega}:=\{x\in\Omega: \ V(x)\neq 0\}.
\eqno(2.18)
$$
Let
$$
\Omega_{n}(a_0,+\infty):=\left\{x\in\Omega: \
U^n\in\mathbb{R}^{m} ,\ a_0\leq|U^n|<+\infty\right\},
$$
where $a_0$ is given in (2.17). Thus, the definition of $\Omega_{n}(a_{0},+\infty)$ and (2.18) imply
$$
\widetilde{\Omega}\subset\Omega_{n}(a_{0},+\infty), \quad \forall n \ \mbox{is large}.
\eqno(2.19)
$$
Note that $(G_1)$, $(G_2)$ and $G\in C^1(\Omega\times\mathbb{R}^m,\mathbb{R})$ imply
$$
|G(x,U)|\leq C_6|U|^2,\quad \forall|U|\leq a_0, \ \ \forall x\in\Omega,
$$
which together with (2.4), (2.8), (2.18), (2.19), $\|V^n\|^2_2\le \tau_2^2$ (below (2.11)), the Fadou's Lemma, $G(x,U)\ge0$ for all $|U|\ge a_0$ in $(2.17)$, $V\neq0$ and $(G_4)$ implies
\\
$$\begin{array}{rcl}
0=\mathlarger{\lim_{n\rightarrow\infty}}\frac{\Phi_{\theta_n}(U^n)}{\|U^n\|^2}
&=&\mathlarger{\lim_{n\rightarrow\infty}}\left[\frac{1}{2}-\mathlarger{\int_{\Omega}}\frac{G(x,U^n)}{\|U^n\|^2}dx
-\frac{\theta\mathlarger{\int_{\Omega}\sum_{j=1}^{m}}h_j(x)u^n_jdx}{\|U^n\|^2}\right]\\
\\
&=&\mathlarger{\lim_{n\rightarrow\infty}}\left[\frac{1}{2}
-\mathlarger{\int_{x\in\Omega_{n}(0,a_{0})}}\frac{G(x,U^n)}{|U^n|^{2}}|V^n|^{2}dx
-\mathlarger{\int_{x\in\Omega_{n}(a_{0},+\infty)}}\frac{G(x,U^n)}{|U^n|^{2}}|V^n|^{2}dx\right]\\
\\
&\leq&\mathlarger{\lim_{n\rightarrow\infty}}\left[\frac{1}{2}+C_6\|V^n\|^{2}_2
-\mathlarger{\int_{x\in\Omega_{n}(a_{0},+\infty)}}\frac{G(x,U^n)}{|U^n|^{2}}|V^n|^{2}dx\right]\\
\\
&\leq&\left(\frac{1}{2}+C_6\tau_2^2\right)
-\mathlarger{\int_{x\in\Omega_{n}(a_{0},+\infty)}\liminf_{n\rightarrow\infty}}
\frac{G(x,U^n)}{|U^n|^{2}}|V^n|^{2}dx\\
\\
&\le&\left(\frac{1}{2}+C_6\tau_2^2\right)
-\mathlarger{\int_{x\in\widetilde{\Omega}}\liminf_{n\rightarrow\infty}}
\frac{G(x,U^n)}{|U^n|^{2}}|V^n|^{2}dx
=-\infty.
\end{array}
\eqno(2.20)
$$
It is a contradiction. Hence, $\{U^n\}$ is bounded in $W$, i.e., $\|U^n\|<\infty$.

\textbf{Part 3.} If the assumptions in Theorem 1.3 hold. If $\|U^n\|\rightarrow\infty$, we let
$$
V^n:=\frac{U^n}{\|U^n\|},
$$
then $\|V^n\|=1$ and (2.5) still holds. Note that (2.4) and (2.7) implies
$$
0\leftarrow\frac{\Phi'_{\theta_n}(U^n)}{\|U^n\|^2}=1- \frac{\mathlarger{\int_{\Omega}}(\nabla G(x,U^n),U^n)dx}{\|U^n\|^2}
-\frac{\theta\mathlarger{\int_{\Omega}\sum_{j=1}^{m}}h_j(x)u^n_jdx}{\|U^n\|^2}.\\
\eqno(2.21)
$$

Suppose $V=0$. Since $(G_1)$ and $(G_7)$ imply $\frac{|\nabla G(x,U)|}{|U|}\leq C_7$ for some $C_7>0$, $\forall (x,U)\in \Omega\times \emph{R}^m$, it follows from the Sobolev compact embedding theorem, (2.5) and $V=0$ that
$$\begin{array}{rcl}
\left|\frac{\mathlarger{\int_{\Omega}(\nabla G(x,U^n),U^n)dx}}{\|U^n\|^2}\right|&=&\left|\frac{\mathlarger{\int_{\Omega}(\nabla G(x,U^n),U^n)|V^n|^2}}{|U^n|^2}dx\right|\\
\\
&\leq&\mathlarger{\int_{\Omega}}\frac{|\nabla G(x,U^n)|}{|U^n|}|V^n|^2dx\\
\\
&\leq&C_7\|V^n\|_2^2\rightarrow0.
\end{array}
\eqno(2.22)
$$
Thus, (2.8), (2.21) and (2.22) imply a contradiction.

Suppose $V\neq0$.
By (2.10), $(G_2)$, $(G_8)$, the Fatou's lemma, the fact $|U_{n}|\rightarrow\infty$ by $V\neq0$ and $\|U_{n}\|\rightarrow\infty$, we have
$$\begin{array}{rcl}
0&=&\mathlarger{\lim_{n\rightarrow\infty}}\frac{\Phi_{\theta_n}(U^n)-\frac{1}{2}\Phi'_{\theta_n}(U^n)U^n}{\|U_n\|^{\mu}}
\\
\\
&\ge&\mathlarger{\lim_{n\rightarrow\infty}}\frac{\mathlarger{\int_{\Omega}}\widetilde{G}(x,U^n)dx-C_1\|U^n\|}{\|U_n\|^{\mu}}\\
\\
&\geq&\mathlarger{\int_{\Omega}\lim_{n\rightarrow\infty}}
\frac{\widetilde{G}(x,U^n)}{|U^n|^{\mu}}|V^n|^{\mu}dx
-\mathlarger{\lim_{n\rightarrow\infty}}\frac{C_1}{\|U^n\|^{\mu-1}}\rightarrow+\infty \qquad (\mu>1),
\end{array}
\eqno(2.23)
$$
which is still a contradiction. Thus, $\{U^n\}$ is bounded in $W$, i.e., $\|U^n\|<\infty$.

(2) Last, we prove $\{U^n\}$ has a convergent subsequence. The boundedness of $\{U^n\}$ in $W$ implies there is a subsequence also
denoted by $\{U^n=(u^n_1,\cdots,u^n_m)\}$ such that
$$
U^n\rightharpoonup U \ \ \mbox{in} \ W, \quad  U^n\to U \ \mbox{a.e.} \ x\in\Omega,
\quad u^n_j\rightharpoonup u_j \ \ \mbox{in} \ H_{0}^{1}(\Omega), \quad u^n_j\to u_j \ \mbox{a.e.} \ x\in\Omega,
\eqno(2.24)
$$
where $j=1,\cdots,m$.
By (2.3) and (2.4), we have
$$\begin{array}{rcl}
0&\leftarrow&\Phi'_{\theta_n}(U^n)(U^n-U)\\
\\
&=&\langle U^n,U^n-U\rangle-\mathlarger{\int_{\Omega}}(\nabla G(x,U^n),U^n-U)dx-\theta\mathlarger{\int_{\Omega}\sum_{j=1}^{m}}h_j(x)(u^n_j-u_j)dx.
\end{array}
\eqno(2.25)
$$
By (2.24), we have
$$
\lim_{n\to\infty}\langle U^n,U^n-U\rangle=\lim_{n\to\infty}\|U^n\|^2-\|U\|^2.
\eqno(2.26)
$$
Note that $(G_{1})$-$(G_{2})$ and $(G_{5})$ or $(G_{7})$ imply that for
$\varepsilon_0$ given in $(G_{1})$ there exists a
$C_{\varepsilon_0}>0$ such that
$$
|\nabla G(x,U^n)|\leq\varepsilon_0|U|+C_{\varepsilon_0}|U|^{p-1}, \quad |G(x,U^n)|\leq\frac{\varepsilon_0}{2}|U|^2+\frac{C_{\varepsilon_0}}{2}|U|^{p}, \quad \forall (x,U)\in \Omega\times \emph{R}^m.
\eqno(2.27)
$$

{\bf Remark of $p$ in (2.27)}: {\it It is fixed by $(G_{5})$ for some $p\in(2,2^*)$ in the super-quadratic case $(G_{4})$, but it can be hold for any $p\in(2,+\infty)$ in the asymptotically-quadratic case $(G_{7})$. Hence, $p$ has a different meaning in the following proofs.}

It follows from (2.27), the H\"{o}lder's inequality, $U^n-U\to 0$ in $[L^2]^m$ and $[L^p]^m$ (by (2.24) and the Sobolev compact embedding theorem), the Sobolev embedding theorem and $\|U^n\|<\infty$ that
$$
\left|\int_{\Omega}(\nabla G(x,U^n),U^n-U)dx\right|\le\varepsilon_0\|U^n\|_2\|U^n-U\|_2+C_{\varepsilon_0}\|U^n\|_p^{p-1}\|U^n-U\|_p\to0.
\eqno(2.28)
$$
By (2.1), the H\"{o}lder's inequality, $U^n-U\to0$ in $L^\mu$ (by (2.24) and the Sobolev compact embedding theorem) and $(H_1)$, we have
$$
\left|\int_{\Omega}\sum_{j=1}^{m}h_j(x)(u^n_j-u_j)dx\right|\le \left(\sum_{j=1,\cdots,m}\|h_j\|_{L^{\mu'}}\right)\|U^n-U\|_{{\mu}}\to0.
\eqno(2.29)
$$
Therefore,  (2.25), (2.26), (2.28) and (2.29) imply $\mathlarger{\lim_{n\to\infty}}\|U^n\|=\|U\|$, which implies $U^n\to U$ in $W$. \ $\Box$\\
\\
\textbf{Lemma 2.3.} {\it $\Phi_{\theta}(U)$ satisfies the conditions $(A2)$-$(A4)$. }\\
\textbf{Proof.} We divide the proof into three parts.

(1) We prove $\Phi_{\theta}(U)$ satisfies the condition $(A2)$. By (2.2), the H\"{o}lder's inequality,  the Sobolev embedding theorem, $(H_1)$ and (2.1), we have
$$
\left|\frac{\partial\Phi_{\theta}(U)}{\partial \theta}\right|=\left|-\int_{\Omega}\sum_{i=1}^{m}h_i(x)u_idx\right|
\le\sum_{i=1}^{m}\|h_i\|_{L^{\mu'}}\|u_i\|_{L^{\mu}}
\le C_{8}\|U\|
$$
for all $(\theta,U)\in[0,1]\times W$. Clearly, $\Phi_{\theta}(U)$ satisfies the condition $(A2)$.

(2) We prove $\Phi_{\theta}(U)$ satisfies the condition $(A3)$.
If $\Phi'_{\theta}(U)=0$, $\forall(\theta,U)\in[0,1]\times W$, then by (2.1)-(2.3), the H\"{o}lder's inequality,  the Sobolev embedding theorem, $\widetilde{G}(x,U)\ge0$ for all $(x,U)\in\Omega\times\mathbb{R}^m$ (by $(G_2)$) and the fact $\widetilde{G}(x,U)\ge C_{9}|U|^\mu$ if $|U|\ge C_{10}$ is large by $(G_3)$ (or $(G_8)$), we have
$$\begin{array}{rcl}
\Phi_{\theta}(U)&=&\Phi_{\theta}(U)-\frac{1}{2}\Phi'_{\theta}(U)U\\
\\
&=&\mathlarger{\int_{\Omega}}\widetilde{G}(x,U)dx
-\frac{\theta}{2}\mathlarger{\int_{\Omega}\sum_{j=1}^{m}}h_j(x)u_jdx\\
\\
&\ge&C_{9}\mathlarger{\int_{\{x\in\Omega: \ |U|\ge C_{10}\}}}|U|^\mu dx-\frac{1}{2}\mathlarger{\sum_{j=1}^{m}}\|h_j\|_{L^{\mu'}}\|u_j\|_{L^\mu}\\
\\
&=&C_{9}\|U\|_\mu^\mu-C_{9}\mathlarger{\int_{\{x\in\Omega: \ |U|\le C_{10}\}}}|U|^\mu dx-\frac{1}{2}\mathlarger{\sum_{j=1}^{m}}\|h_j\|_{L^{\mu'}}\|u_j\|_{L^\mu}\\
\\
&\ge&C_{9}\|U\|_\mu^\mu-C_{9}C_{10}^\mu meas(\Omega)-\frac{1}{2}\left(\mathlarger{\sum_{j=1,\cdots,m}}\|h_j\|_{L^{\mu'}}\right)\|U\|_\mu\\
\\
&\ge&\frac{C_{9}}{2}\|U\|_\mu^\mu-C_{11}, \quad \forall(\theta,U)\in[0,1]\times W.
\end{array}$$
That is,
$$\begin{array}{rcl}
\|U\|_\mu^{2\mu}\le \left(\frac{2}{C_{9}}\right)^2 \left(\Phi_{\theta}(U)+C_{11}\right)^2
&\le& 2\left(\frac{2}{C_{9}}\right)^2 \left(\Phi_{\theta}(U)^2+C_{11}^2\right)\\
\\
&\le& C_{12}\left(\Phi_{\theta}(U)^2+1\right), \quad \forall(\theta,U)\in[0,1]\times W.
\end{array}$$
It follows from the H\"{o}lder's inequality, $(H_1)$ and (2.1) that
$$\begin{array}{rcl}
\left|\frac{\partial\Phi_{\theta}(U)}{\partial \theta}\right|&\le&\mathlarger{\sum_{j=1}^{m}}\|h_j\|_{L^{\mu'}}\|u_j\|_{L^{\mu}}\\
\\
&\le& \left(\mathlarger{\sum_{j=1,\cdots,m}}\|h_j\|_{L^{\mu'}}\right)\|U\|_\mu\\
\\
&\le& C_{13}\left(\Phi^2_{\theta}(U)+1\right)^{\frac{1}{2\mu}}.
\end{array}$$
Therefore, $\Phi_{\theta}$ satisfies the condition $(A3)$ with
$$
\eta_1(\theta,s)=-\eta_2(\theta,s)=-C_{13}\left(s^2+1\right)^{\frac{1}{2\mu}}.
$$

(3) We prove $\Phi_{\theta}(U)$ satisfies the condition $(A4)$. Obviously, $\Phi_{0}(-U)=\Phi_{0}(U)$.
For $U\in H\subset W$ with $\dim H<\infty$, there exists $\widehat{\zeta}>0$ such that
$$
\|U\|^2\le 2\widehat{\zeta}\|U\|^2_2.
$$
Thus, from $(G_4)$ or $(G_7)$, we can take $\zeta_0=2\widehat{\zeta}$ in (2.17). It follows from (2.1), (2.2), (2.17), $(H_1)$, the H\"{o}lder's inequality and $\|U\|^2\le 2\widehat{\zeta}\|U\|^2_2$, we have
$$\begin{array}{rcl}
\Phi_{\theta}(U)
&=&\frac{1}{2}\|U\|^2-\mathlarger{\int_{\Omega}}G(x,U)dx
-\theta\mathlarger{\int_{\Omega}\sum_{j=1}^{m}}h_j(x)u_jdx\\
\\
&=&\frac{1}{2}\|U\|^2-\left(\mathlarger{\int_{\{x\in\Omega: \ |U(x)|\ge a_0\}}}+\mathlarger{\int_{\{x\in\Omega: \ |U(x)|< a_0\}}}\right)G(x,U)dx
-\theta\mathlarger{\int_{\Omega}\sum_{j=1}^{m}}h_j(x)u_jdx\\
\\
&\le&\widehat{\zeta}\|U\|^2_2-2\widehat{\zeta}\mathlarger{\int_{\{x\in\Omega: \ |U(x)|\ge a_0\}}}|U|^2dx+meas(\Omega)\mathlarger{\max_{x\in\Omega,|U|\le a_0}}|G(x,U)|\\
\\
&&+\theta \|U\|_\mu\mathlarger{\sum_{j=1,\cdots,m}}\|h_j\|_{L^{\mu'}}\\
\\
&=&\widehat{\zeta}\|U\|^2_2-2\widehat{\zeta}\mathlarger{\int_{\Omega}}|U|^2dx
+2\widehat{\zeta}\mathlarger{\int_{\{x\in\Omega: \ |U(x)|< a_0\}}}|U|^2dx+meas(\Omega)\mathlarger{\max_{x\in\Omega,|U|\le a_0}}|G(x,U)|\\
\\
&&+\theta \|U\|_\mu\mathlarger{\sum_{j=1,\cdots,m}}\|h_j\|_{L^{\mu'}}\\
\\
&\le&-\widehat{\zeta}\|U\|^2_2+\underbrace{2\widehat{\zeta}a_0^2 meas(\Omega)
+meas(\Omega)\max_{x\in\Omega,|U|\le a_0}|G(x,U)|}
+\theta \|U\|_\mu\mathlarger{\sum_{j=1,\cdots,m}}\|h_j\|_{L^{\mu'}}\\
\\
&=&-\widehat{\zeta}\|U\|^2_2+\underbrace{\widetilde{C}}
+\left(\theta \mathlarger{\sum_{j=1,\cdots,m}}\|h_j\|_{L^{\mu'}}\right)\|U\|_\mu,\quad \forall U\in H \ \mbox{with} \ \|U\|\to\infty.
\end{array}
$$
It follows from the equivalence of all norms in the finite dimensional
subspace $H$ that
$$
\lim_{U\in H,\ \|U\|\to\infty}\sup_{\theta\in[0,1]}\Phi_{\theta}(U)=-\infty, \quad \forall H\subset W \ \mbox{with} \ \dim H<\infty.
$$
Therefore, $\Phi_{\theta}$ satisfies the condition $(A4)$. \ $\Box$\\

To apply Lemma 2.1, we let $\{S_n\}_{n\ge1}$ be an orthonormal basis of $W$.
Let
$$
W_k:=\mbox{span}\{S_1,S_2,\cdots,S_k\},
$$
$$
W^\bot_{k-1}:=\overline{\mbox{span}\{S_k,S_{k+1},\cdots\}},
$$
and
$$
c_k=\inf_{\gamma\in \Gamma_k}\sup_{U\in W_k}\Phi_0(\gamma(U)), \quad k\in\mathbb{N}^+,
\eqno(2.30)
$$
where
$$
\begin{array}{cl}
\Gamma_k:=&\big\{\gamma\in C(W_k,W):\ \gamma(U)=(\gamma_1(u_1),\cdots,\gamma_m(u_m))=-\gamma(-U), \\
\\
 &   \ \gamma_j(u)=u\in H_{0}^{1}(\Omega) \ for \ \|u\|_{\overline{H}_j} \ large, \ j=1,\cdots,m, \ U=(u_1,\cdots,u_m)\big\}.
\end{array}$$
By the definition of $c_k$, we have
$$
0<c_1\le c_2\le\cdots \le c_k\le\cdots.
\eqno(2.31)
$$

Now, by Lemmas 2.2 and 2.3, we have our functionals $\Phi_{\theta}$ ($\theta\in[0,1]$) defined in (2.2) satisfy the conditions $(A1)$-$(A4)$ of Lemma 2.1 with
$$
\eta_1(\theta,s)=-\eta_2(\theta,s):=-C_{13}\left(s^2+1\right)^{\frac{1}{2\mu}},
\eqno(2.32)
$$
and
$$
\overline{\eta}_1(s)=\overline{\eta}_2(s)=C_{13}\left(s^2+1\right)^{\frac{1}{2\mu}}.
\eqno(2.33)
$$
Moreover, by $\eta_2(\theta,s)\ge0$ (by (2.32)) and $(2.0)$ in Lemma 2.1, we have $\psi_2(\cdot,s)$ is non-decreasing on $[0, 1]$. It follows from $(i)$ of Lemma 2.1 and $(2.0)$ that
$$
c_k =\psi_2(0,c_k)\le\psi_2(1,c_k)<\psi_1(1,c_{k+1})\le \widetilde{c_k},\quad \forall k\in\mathbb{N}^+.
\eqno(2.34)
$$\\
\textbf{Remark 2.1.} Note that Lemma 2.1 have two results: $(i)$ and $(ii)$, but the result $(ii)$ of Lemma 2.1 does not hold for $k$ large enough in our conditions, which will be proved in Section 3. Moreover, we will prove $c_k \to +\infty$ (see (3.12) and (3.14)), then $\widetilde{c_k} \to +\infty$. Therefore, our Theorems 1.1-1.3 can be proved by $(i)$ of Lemma 2.1 for $k$ large enough.\\

The following two lemmas are essential in the proof that the result $(ii)$ of Lemma 2.1 does not hold for $k$ large enough.\\
\\
\textbf{Lemma 2.4} (\cite{Cwickel,Li}). {\it Let $V\in L^{\frac{N}{2}}(\Omega)$, for the eigenvalue problem
$$\left\{
\begin{array}{ll}
-\triangle u-V(x) u=\lambda u,& x\in\Omega\subset\mathbb{R}^{N},\\
u=0,& x\in\partial\Omega,
\end{array}\right.
\eqno(2.35)
$$
we have two results.\\
$(1)$ If $N=2$, then for $\forall\varepsilon>0$ there is $\overline{C_\varepsilon} > 0$ such that
$$
\sharp_-(-\triangle -V)\le \overline{C_\varepsilon}\|V\|_{L^{\varepsilon+1}}^{\varepsilon+1},
$$
where $\sharp_-(-\triangle -V)$ denotes the number of non-positive eigenvalues of $(2.31)$.\\
$(2)$ If $N\geq3$, then there is $C_N > 0$ such that }
$\sharp_-(-\triangle -V)\le C_N\|V\|_{L^{\frac{N}{2}}}^{\frac{N}{2}}.$\\
\\
\textbf{Lemma 2.5} (\cite{Lazer,Tanaka}). {\it Let $E$ be a Hilbert space with the norm $\|\cdot\|_E$, and $K\in C^2(E, \mathbb{R})$ be a functional satisfies the usual $(PS)$ condition and the following conditions hold.\\
$(K1)$ $K(0)=0$ and $K(-u)=K(u)$, $\forall u\in E$.\\
$(K2)$ For any $E_k\subset E$ with $\dim E_k =k<\infty$, there is $R_k> 0$ such that $K(u)<0$, $\forall u\in E_k$ with $\|u\|_E\ge R_k.$\\
$(K3)$ $K'(u)=u+\kappa(u)$, where $\kappa$ is a compact operator.
\\
Then there exists $\{u_n\} \subset E$ such that
$K(u_k)\le d_k$, $K'(u_k)=0$, $index(K''(u_k))\ge k$,
where }
$$
\mbox{index}(K''(u))=\max\left\{\dim S: \ S\subset E \ such \ that \ K''(u)(l,l)\le0, \ l\in S\right\},
\quad
d_k=\inf_{\gamma\in \widehat{\Gamma}_k}\sup_{u\in D_k}K(\gamma(u)),
$$
$$
D_k:=\left\{u\in E_k: \ \|u\|_E\le R_k\right\}, \quad
\widehat{\Gamma}_k=\left\{\gamma\in C(D_k,E):\ \gamma \ is \ odd, \ \gamma|_{\partial D_k}=id\right\}.
$$

\section * {3. Proofs of the main results}

By Remark 2.1, we need to prove the result $(ii)$ of Lemma 2.1 does not hold for $k$ large enough. To do this, we should give the upper bound and lower bound of the minimax value sequence $\{c_k\}$, which are based on Lemma 5.3 in \cite{Bahri} and the above Lemmas 2.4 and 2.5.\\
\\
\textbf{Proofs of Theorems 1.1-1.3.}
By contradiction, suppose the result $(ii)$ of Lemma 2.1 holds for $k$ large enough, i.e.,
$$
c_{k+1}-c_{k}\le C\left[\overline{\eta}_1(c_{k+1})+\overline{\eta}_2(c_{k})+1\right], \quad \forall k\ge \widetilde{k_0}>0.
$$
Then, it follows from (2.31), (2.33) and $(s+t)^q\le s^q+t^q$ for $s,t\ge0$ and $0<q<1$ that
$$
c_{k+1}-c_k\le C_{14}\left(c_{k+1}^{\frac{1}{\mu}}+c_{k}^{\frac{1}{\mu}}+1\right), \quad \forall k\ge \widetilde{k_0}.
\eqno(3.1)
$$

\textbf{Part 1.} We give the upper bound of the minimax value sequence $\{c_k\}$, i.e., the following claim 1.

\textbf{Claim 1.} We claim
$$
c_k\le C_{15}k^\frac{\mu}{\mu-1}, \quad \forall k\ge k_0>0.
\eqno(3.2)
$$
Next, we prove Claim 1. By the similar proof in Lemma 5.3 of \cite{Bahri}, we let
$$
s:=\frac{1}{\mu}\in(0,1),\qquad \chi_k:=c_k k^{-\frac{1}{1-s}},
$$
then we only need to prove $\chi_k\le D$ for $k$ large enough. In the proof of Claim 1, we let $D$ denote different positive constant.
Note that (2.31) and (3.1) implies
$$
c_{k+1}-c_k\le D c_{k+1}^{s},  \quad \forall k\ge k_1 \ \mbox{for some} \ k_1>0,
$$
it follows from $\chi_k:=c_k k^{-\frac{1}{1-s}}$ that
$$
\begin{array}{rcl}
\chi_{k+1}(k+1)^{\frac{1}{1-s}}-\chi_k k^{\frac{1}{1-s}}&\le& D \chi_{k+1}^{s} (k+1)^{\frac{s}{1-s}}\\
\\
&=&D \chi_{k+1}^{s} (k+1)^{-1}(k+1)^{\frac{1}{1-s}},  \quad \forall k\ge k_1.
\end{array}
$$
Thus, we have
$$\begin{array}{rcl}
\chi_{k+1}(1+\frac{1}{k})^{\frac{1}{1-s}}-\chi_k &\le& D \chi_{k+1}^{s} (k+1)^{-1}(1+\frac{1}{k})^{\frac{1}{1-s}}\\
\\
&\le& D \chi_{k+1}^{s} (k+1)^{-1},  \quad \forall k\ge k_1,
\end{array}$$
which together with the fact $(1+a)^{b}\ge 1+ab$ with $a>0$ and $b=\frac{1}{1-s}>1$ implies
$$
\chi_{k+1}+\chi_{k+1}k^{-1}\frac{1}{1-s}-\chi_k \le D \chi_{k+1}^{s} (k+1)^{-1},  \quad \forall k\ge k_1.
\eqno(3.3)
$$
Suppose that $\{\chi_k\}$ is unbounded, then there exists a subsequence still denoted by $\{\chi_k\}$ such that
$$
\chi_{k+1}\le \chi_k\quad \mbox{or} \quad \chi_{k+1}\ge \chi_k,\quad \forall k\ge k_1.
$$
If $\chi_{k+1}\le \chi_k$, $\forall k\ge k_1$, then $\chi_k\le \chi_{k_1}$, $\forall k\ge k_1$.
If $\chi_{k+1}\ge \chi_k$, $\forall k\ge k_1$, then it follows from (3.3) that
$$
\chi_{k+1}k^{-1}\frac{1}{1-s} \le D \chi_{k+1}^{s} (k+1)^{-1},  \quad \forall k\ge k_1,
$$
thus $\chi_k\le\chi_{k+1}\le D \chi_{k+1}^{s}$, i.e., $\chi_k\le\chi_{k+1}\le D_0$ for some $D_0>0$.
Hence,
$$
\chi_{k}\le \max\{D_0,\chi_{k_1}\},\quad \forall k\ge k_1>0,
$$
which contradicts with $\{\chi_k\}$ is unbounded. Therefore, $\chi_k\le C$ for $k$ large enough. That is, Claim 1 holds.

\textbf{Part 2.} We give the lower bound of the minimax value sequence $\{c_k\}$ and give the contradiction when $(ii)$ of Lemma 2.1 holds.

Since (2.27) and the definition of $\|U\|$ imply
$$\begin{array}{rcl}
\Phi_0(U)&=&\frac{1}{2}\|U\|^2-\mathlarger{\int_{\Omega}}G(x,U)dx\\
\\
&\ge& \frac{1}{2}\|U\|^2- \frac{\varepsilon_0}{2}\|U\|_2^2-\frac{C_{\varepsilon_0}}{2}\|U\|_p^p\\
\\
&=& \frac{1}{2}\left(\|u_1\|_{\overline{H}_1}^2+\cdots+\|u_m\|_{\overline{H}_m}^2\right)
- \frac{\varepsilon_0}{2}\mathlarger{\int_{\Omega}}(u_1^2+\cdots+u_m^2)dx
-\frac{C_{\varepsilon_0}}{2}\mathlarger{\int_{\Omega}}(u_1^2+\cdots+u_m^2)^{\frac{p}{2}}dx\\
\\
&:=& \widehat{K}(U)=\widehat{K}(u_1,\cdots,u_m), \quad \forall U=(u_1,\cdots,u_m)\in W,
\end{array}
\eqno(3.4)
$$
which implies
$$\begin{array}{rcl}
c_k&=&\mathlarger{\inf_{\gamma\in \Gamma_k}\sup_{U\in W_k}}\Phi_0(\gamma(U))\\
\\
&\ge&\mathlarger{\inf_{\gamma\in \Gamma_k}\sup_{U\in W_k}}\widehat{K}(\gamma(U))\\
\\
&\ge&\mathlarger{\inf_{\gamma\in \Gamma_k}\sup_{U=(0,\cdots,u_j,\cdots,0)\in W_k}}\widehat{K}(0,\cdots,\gamma_j(u_j),\cdots,0), \quad j=1,\cdots,m.
\end{array}
\eqno(3.5)
$$

Let
$$
K_j(u):=\frac{1}{2}\|u\|_{\overline{H}_j}^2- \frac{\varepsilon_0}{2}\int_{\Omega}|u|^2dx-\frac{C_{\varepsilon_0}}{2}\int_{\Omega}|u|^pdx, \quad u \in H_0^1(\Omega), \quad j=1,\cdots,m,
\eqno(3.6)
$$
then
$$
K_1(u_1)=\widehat{K}(u_1,0,\cdots,0), \quad K_2(u_2)=\widehat{K}(0,u_2,0,\cdots,0),  \quad \cdots , \quad K_m(u_m)=\widehat{K}(0,\cdots,0,u_m).
\eqno(3.7)
$$
Hence, $K_j(u)$ ($j=1,\cdots,m$) satisfies $(K1)$-$(K3)$ in Lemma 2.5 with $E=H_0^1(\Omega)$, $[D_k]^m \subset [E_k]^m =W_k$ and $R_k$ large enough, which together with (3.5) and (3.7) implies
$$\begin{array}{rcl}
mc_k&\ge&\mathlarger{\inf_{\gamma\in \Gamma_k}\sup_{U=(u_1,0,\cdots0)\in W_k}}\widehat{K}(\gamma_1(u_1),\cdots,0)\\
\\
&\quad&+\cdots+\mathlarger{\inf_{\gamma\in \Gamma_k}\sup_{U=(0,\cdots,u_m)\in W_k}}\widehat{K}(0,\cdots,\gamma_m(u_m))\\
\\
&=&\mathlarger{\inf_{\gamma\in \Gamma_k}\sup_{u_1\in E_k}}K_1(\gamma_1(u_1)) +\cdots
+\mathlarger{\inf_{\gamma\in \Gamma_k}\sup_{u_m\in E_k}}K_m(\gamma_m(u_m))\\
\\
&\ge&\mathlarger{\inf_{\gamma_1\in \widehat{\Gamma}_k}\sup_{u_1\in D_k}K_1(\gamma_1(u_1))+\cdots+\inf_{\gamma_m\in \widehat{\Gamma}_k}\sup_{u_m\in D_k}K_m(\gamma_m(u_m))}\\
\\
&\ge&md_k.
\end{array}
\eqno(3.8)
$$
By applying Lemma 2.5 to $K_j$ with the corresponding $d_k$, we have that there exists
$\{u_k\} \subset H_0^1(\Omega)$ such that
$$
K_j(u_k)\le d_k,
\eqno(3.9)
$$
$$
K'_j(u_k)=0,
\eqno(3.10)
$$
and
$$
\mbox{index}(K''_j(u_k))\ge k, \quad  j=1,\cdots,m.
\eqno(3.11)
$$
The definition of $K_j$ implies
$$\begin{array}{rcl}
K''_j(u_k)(l,l)&=&\left(\left[-\triangle+(V_j(x)-\varepsilon_0) I-\frac{C_{\varepsilon_0}}{2}p(p-1)|u_k|^{p-2}\right]l,l\right), \\
\\
&&\forall u_k\in H_0^1(\Omega), \quad  j=1,\cdots,m,
\end{array}
\eqno(3.12)
$$
where $(\cdot,\cdot)$ denotes the duality product between $H^{-1}(\Omega)$ and $H_0^1(\Omega)$. Thus, by (3.11), (3.12) and $\varepsilon_0<\underline{V}\le V_j(x)$ (by $(V_1)$ and $(G_1)$), we have
$$\begin{array}{rcl}
k&\le&\sharp_-\left[-\triangle+(V_j(x)-\varepsilon_0) I-\frac{C_{\varepsilon_0}}{2}p(p-1)|u_k|^{p-2}\right]\\
\\
&\le&\sharp_-\left[-\triangle-\frac{C_{\varepsilon_0}}{2}p(p-1)|u_k|^{p-2}\right].
\end{array}
\eqno(3.13)
$$
By (3.6) and (3.8)-(3.10), we have
$$
c_k\ge d_k\ge K_j(u_k)-\frac{1}{2}K'_j(u_k)u_k= \frac{p-2}{4}C_{\varepsilon_0}\|u_k\|_{L^{p}}^p.
\eqno(3.14)
$$

\textbf{Case 1.} If $N=2$. By (3.13) and (1) in Lemma 2.4 with
$$
V(x)=\frac{1}{2}C_{\varepsilon_0}p(p-1)|u_k|^{p-2},
$$
we have
$$
k\le \overline{C_\varepsilon}\|\frac{1}{2}C_{\varepsilon_0}p(p-1)|u_k|^{p-2}\|_{L^{\varepsilon+1}}^{\varepsilon+1}
=C_{16}\|u_k\|_{L^{(p-2)(\varepsilon+1)}}^{(p-2)(\varepsilon+1)},
\eqno(3.15)
$$
where $C_{16}:=\overline{C_\varepsilon}\left[\frac{1}{2}C_{\varepsilon_0}p(p-1)\right]^{\varepsilon+1}$ and $(p-2)(\varepsilon+1)\ge1$, i.e.,
$$
p\ge 2+\frac{1}{\varepsilon+1}.
\eqno(3.16)
$$
Suppose that $p\ge(p-2)(\varepsilon+1)$, i.e.,
$$
p\le\frac{2(\varepsilon+1)}{\varepsilon},
\eqno(3.17)
$$
then we have $L^p(\Omega) \subset L^{(p-2)(\varepsilon+1)}(\Omega)$ and
$$
\|u_k\|_{L^{(p-2)(\varepsilon+1)}}\le C_{17}\|u_k\|_{L^{p}}.
\eqno(3.18)
$$
By (3.14), (3.15) and (3.18), we have
$$
k\le C_{16}\|u_k\|_{L^{(p-2)(\varepsilon+1)}}^{(p-2)(\varepsilon+1)}\le C_{16}C_{17}^{(p-2)(\varepsilon+1)}\|u_k\|_{L^{p}}^{(p-2)(\varepsilon+1)}\le C_{18}c_k^{\frac{(p-2)(\varepsilon+1)}{p}},
$$
that is,
$$
c_k\ge C_{19} k^{\frac{p}{(p-2)(\varepsilon+1)}}.
\eqno(3.19)
$$

$(i)$ Super-quadratic case.
By (3.2) and (3.19), we can get a contradiction if $\frac{p}{(p-2)(\varepsilon+1)}>\frac{\mu}{\mu-1}$ for $k$ is large enough, i.e.,
$$
p<\frac{2\mu(\varepsilon+1)}{\mu\varepsilon+1}.
\eqno(3.20)
$$
Note that
$$
\frac{2\mu(\varepsilon+1)}{\mu\varepsilon+1}<\frac{2(\varepsilon+1)}{\varepsilon},
$$
thus it follows from (3.16), (3.17) and (3.20) that we can get a contradiction for $k$ is large enough if and only if
$$
2+\frac{1}{\varepsilon+1}<\frac{2\mu(\varepsilon+1)}{\mu\varepsilon+1} \quad \mbox{and}\quad
2+\frac{1}{\varepsilon+1}\le p<\frac{2\mu(\varepsilon+1)}{\mu\varepsilon+1}, \qquad \forall\varepsilon>0,
$$
that is, we get a contradiction for $k$ is large enough if and only if
$$
\mu>1+\frac{\varepsilon+1}{\varepsilon+2}\quad \mbox{and}\quad
2+\frac{1}{\varepsilon+1}\le p<\frac{2\mu(\varepsilon+1)}{\mu\varepsilon+1}, \qquad \forall\varepsilon>0,
\eqno(3.21)
$$
which the condition (1.5) in (1) of Theorem 1.1.

$(ii)$ Asymptotically-quadratic case. Note that we can take
$$
\forall p\in\left[2+\frac{1}{\varepsilon+1},\frac{2\mu(\varepsilon+1)}{\mu\varepsilon+1}\right)
$$
in the asymptotically-quadratic case such that the second inequality holds in (3.21), so there is a contradiction for $k$ is large enough by (3.21) if and only if
$$
\mu>1+\frac{\varepsilon+1}{\varepsilon+2}, \qquad \forall\varepsilon>0,
$$
which the condition (1.7) in (1) of Theorem 1.3.

\textbf{Case 2.} If $N\ge3$. By (3.13) and (2) in Lemma 2.4 with
$$
V(x)=\frac{1}{2}C_{\varepsilon_0}p(p-1)|u_k|^{p-2},
$$
we have
$$
k\le C_N\|\frac{1}{2}C_{\varepsilon_0}p(p-1)|u_k|^{p-2}\|_{L^{\frac{N}{2}}}^{\frac{N}{2}}
=C_{20}\|u_k\|_{L^{\frac{(p-2)N}{2}}}^{\frac{(p-2)N}{2}},
\eqno(3.22)
$$
where $C_{20}:=C_N\left[\frac{1}{2}C_{\varepsilon_0}p(p-1)\right]^{\frac{N}{2}}$ and $\frac{(p-2)N}{2}\ge1$, i.e.,
$$
p\ge2+\frac{2}{N}.
\eqno(3.23)
$$
Since $p<2^\ast$ implies $\frac{(p-2)N}{2}< p$,
we have $L^p(\Omega) \subset L^{\frac{(p-2)N}{2}}(\Omega)$ and
$$
\|u_k\|_{L^{\frac{(p-2)N}{2}}}\le C_{21}\|u_k\|_{L^{p}}.
\eqno(3.24)
$$
By (3.14), (3.22) and (3.24), we have
$$
k\le C_{20}\|u_k\|_{L^{\frac{(p-2)N}{2}}}^{\frac{(p-2)N}{2}} \le C_{20}C_{21}^{\frac{(p-2)N}{2}}\|u_k\|_{L^{p}}^{\frac{(p-2)N}{2}}\le C_{22} c_k^{\frac{(p-2)N}{2p}},
$$
that is,
$$
c_k\ge C_{23} k^{\frac{2p}{(p-2)N}}.
\eqno(3.25)
$$

$(i)$ Super-quadratic case. By (3.2) and (3.25), we can get a contradiction if $\frac{2p}{(p-2)N}>\frac{\mu}{\mu-1}$ for $k$ is large enough, i.e.,
$$
p<\frac{2\mu N}{\mu(N-2)+2}.
$$
It follows from (3.23) that we can get a contradiction for $k$ is large enough if and only if
$$
2+\frac{2}{N}<\frac{2\mu N}{\mu(N-2)+2}\quad \mbox{and}\quad
2+\frac{2}{N}\le p<\frac{2\mu N}{\mu(N-2)+2}.
$$
That is, we get a contradiction for $k$ is large enough if and only if
$$
\mu>1+\frac{N}{N+2}\quad \mbox{and}\quad 2+\frac{2}{N}\le p<\frac{2\mu N}{\mu(N-2)+2},
\eqno(3.26)
$$
which the condition (1.6) in (2) of Theorem 1.1.

$(ii)$ Asymptotically-quadratic case. Note that we can take
$$
\forall p\in\left[2+\frac{2}{N},\frac{2\mu N}{\mu(N-2)+2}\right)
$$
in the asymptotically-quadratic case such that the second inequality holds in (3.26), so we get a contradiction for $k$ is large enough by (3.26) if and only if
$$
\mu>1+\frac{N}{N+2},
$$
which the condition (1.8) in (2) of Theorem 1.3.

Therefore, the proofs of Theorems 1.1-1.3 are completed by Remark 2.1. \ $\Box$

\section*{Conflicts of interest}

The author declares that he has no conflict of interest.

\section*{Data availability statement}
My manuscript has no associated date.

\end{document}